\documentclass[12pt]{amsart}

\usepackage{amsxtra,amscd}
\usepackage[dvips]{graphics}

\usepackage{geometry}
\usepackage{amsfonts}
\usepackage{amssymb}
\usepackage[mathscr]{eucal}
\usepackage{amsmath}
\usepackage{amsthm}
\usepackage[all]{xy}
\usepackage[mathscr]{eucal}
\usepackage{graphicx}
\usepackage{color}
\definecolor{red-}{rgb}{1.0,0.0,0.0}
\definecolor{green-}{rgb}{0.0,0.7,0.0}
\definecolor{brown-}{rgb}{0.9,0.6,0.0}

\usepackage{hyperref}
\hypersetup{
    colorlinks=true,
    linkcolor=blue,
    citecolor=blue,
    filecolor=blue,
    urlcolor=blue
}

\newtheorem{defi}{Definition}[section]

\newtheorem{prop}[defi]{Proposition}

\newtheorem{quest}[defi]{Question}

\begin{document}

\title[The Hilbert curve of a $4$-dimensional scroll]{The Hilbert curve of a $4$-dimensional scroll\\
 with a divisorial fiber}

\author{Antonio Lanteri and Andrea Luigi Tironi}

\date{\today}

\address{Dipartimento di Matematica ``F. Enriques'',
Universit\`a degli Studi di Milano, Via C. Saldini, 50,  I-20133 Milano,
Italy} \email{antonio.lanteri@unimi.it}

\address{
Departamento de Matem\'atica, Universidad de Concepci\'on, Casilla
160-C, Concepci\'on, Chile} \email{atironi@udec.cl}


\begin{abstract}
In dimension $n = 2m-2 \geq 4$ adjunction theoretic scrolls over a smooth $m$-fold may not be
classical scrolls, due to the existence of divisorial fibers. A $4$-dimensional scroll $(X,L)$ over $\mathbb P^3$ of this type is considered, and the equation of its Hilbert curve $\Gamma$ is determined in two ways, one of which relies on the fact that $(X,L)$ is at the same time a classical scroll over a threefold $Y \not=\mathbb P^3$. It turns out that $\Gamma$ does not perceive divisorial fibers. The equation we obtain also shows that a question raised in \cite{BLS} has
negative answer in general for non-classical scrolls over a $3$-fold. More precisely, the answer for
$(X,L)$ is negative or positive according to whether $(X,L)$ is regarded as an adjunction theoretic scroll or as a classical scroll; in other words, it is the
answer to this question to distinguish between the existence of jumping fibers or not.
\end{abstract}

\maketitle

\section*{Introduction}\label{Intro}

For a polarized manifold $(\mathcal X, \mathcal L)$ of dimension $n$, two notions of scroll over a variety $\mathcal Y$ of
smaller dimension $m$ are possible: $(\mathcal X,\mathcal L)$ is a {\it{classical scroll}} if $\mathcal X=\mathbb P(\mathcal E)$ for an ample vector bundle $\mathcal E$ on $\mathcal Y$, $\mathcal L$ being
the tautological line bundle, while $(\mathcal X,\mathcal L)$ is an {\it{adjunction theoretic scroll}} over $\mathcal Y$ if there exists a surjective
morphism $\varphi:\mathcal X \to \mathcal Y$ such that $K_{\mathcal X} + (n-m+1)\mathcal L = \varphi^*\mathcal A$ for some ample line bundle
$\mathcal A$ on $\mathcal Y$ (see \cite[p.\ 81]{BS}).
Essentially, classical scrolls are also adjunction theoretic scrolls, by taking as $\varphi$ the bundle projection
$p:\mathcal X \to \mathcal Y$, except when $K_{\mathcal Y}+\det \mathcal E$ fails to be ample, and all these exceptions are well known in low dimension (see, \cite{BS0}, \cite[$\S 3$]{T}
and \cite[$\S 4.2$]{T1}). Conversely, it is known that for $m \leq 4$ an adjunction theoretic scroll is a classical
scroll if $n \geq 2m-1$, when $\mathcal L$ is very ample (see \cite[Proposition 14.1.3]{BS} and \cite[Theorem 2.2]{T}). This is no longer
true for $n = 2m-2 \geq 4$, since in this case $\varphi$ can admit divisorial fibers. A class of examples illustrating this phenomenon
is due to Beltrametti and Sommese \cite[(4.2)]{BS0}.

In this paper, a $4$-dimensional scroll $(X,L)$ over $\mathbb P^3$ -- the simplest example
of this type -- is considered and the equation of its Hilbert curve is determined.
This is done in two different ways: the former is via the explicit Riemann--Roch formula
for $4$-folds exploiting that $(X,L)$ itself is also a classical scroll over another threefold $Y$, related to $\mathbb P^3$
(Section \ref{first});
the latter relies on a recursive procedure introduced in \cite[Section 4]{LT}, working for scrolls of both types (Section \ref{second}). It turns out that the Hilbert curve does not detect
divisorial fibers. Moreover, the equation we obtain indicates that a question raised in \cite{BLS} has negative answer in general
for non-classical scrolls. More precisely, it turns out that for our $(X,L)$ the answer is negative or positive according to whether
we look at it either as an adjunction theoretic scroll over $\mathbb P^3$, or as a classical scroll over $Y$; in other words, it is the
answer to this question to distinguish between the existence of jumping fibers or not.

\section{Preliminaries}\label{prel}

Varieties considered in this paper are defined over the field $\mathbb C$ of complex numbers.
We use the standard notation and terminology from algebraic geometry. A manifold is any smooth projective variety.
Tensor products of line bundles are denoted additively. The pullback of a vector bundle $\mathcal E$ on a manifold $\mathcal X$
by an embedding $\mathcal Z \to \mathcal X$ is simply denoted by $\mathcal E_{\mathcal Z}$, while $K_{\mathcal X}$
will stand for the canonical
bundle of $\mathcal X$. A polarized manifold is a pair $(\mathcal X,\mathcal L)$ consisting of a manifold
$\mathcal X$ and an ample line bundle $\mathcal L$ on $\mathcal X$.

\medskip

For the notion and the general properties of the Hilbert curve associated to a
polarized manifold we refer to \cite{BLS}, see also \cite{L}. Here we just recall some basic facts.
Let $(\mathcal X,\mathcal L)$ be a polarized manifold of dimension $n \geq 2$ and regard
$\text{N}(\mathcal X):= \text{Num}(\mathcal X) \otimes_{\mathbb Z} \mathbb C$ as a complex affine space.
If $\text{rk}\langle K_{\mathcal X}, \mathcal L \rangle = 2$, we can consider the plane
$\mathbb A^2= \mathbb C \langle K_{\mathcal X},\mathcal L \rangle \subset \text{N}(\mathcal X)$,
generated by the classes of $K_{\mathcal X}$ and $\mathcal L$.
For any line bundle $D$ on $\mathcal X$ the Riemann--Roch theorem provides an expression for the
Euler--Poincar\'e  characteristic $\chi(D)$ in terms of $D$ and the Chern classes of $\mathcal X$.
Let $p$ denote the complexified polynomial of $\chi(D)$, when we set $D = xK_{\mathcal X} + y\mathcal L$,
with $x,y$ complex numbers, namely $p(x,y) = \chi (xK_{\mathcal X}+y\mathcal L)$.
The Hilbert curve of $(\mathcal X,\mathcal L)$ is the complex affine plane curve $\Gamma = \Gamma_{(\mathcal X,\mathcal L)} \subset \mathbb A^2$ of
degree $n$ defined by $p(x,y)=0$ \cite[Section 2]{BLS}.
Notice that the Hilbert curve can be defined also when the numerical classes of $K_{\mathcal X}$ and $\mathcal L$ are linearly
dependent, but in this case, the $(x,y)$-plane is only formal, $\Gamma_{(\mathcal X,\mathcal L)}$ losing the meaning of
a plane section of the Hilbert variety of $\mathcal X$ (see \cite[Section 2]{BLS}).
For example, the Hilbert curve of $\big(\mathbb P^n, \mathcal O_{\mathbb P^n}(r)\big)$ has the following equation (see, e.\ g., \cite[p.\ 465]{BLS} and \cite[Theorem 2.7]{LT}):
\begin{equation}\label{p1}
p(x,y) = \frac{(-1)^n}{n!} \prod_{i=1}^n \big((n+1)x-ry-i \big)\ .
\end{equation}

Due to Serre duality, $\Gamma$ is invariant under the involution $D \mapsto K_{\mathcal X}-D$ acting on $\text{N}(\mathcal X)$.
Sometimes, to make this symmetry more evident, it is convenient to represent $\Gamma$ in terms of the affine coordinates $(u=x - \frac{1}{2},
v =y)$ rather than $(x,y)$. So, rewriting our divisor as $D=\frac{1}{2}K_{\mathcal X} + \Delta$, where $\Delta=uK_{\mathcal X}+v\mathcal L$,
$\Gamma$ can be represented with respect to these coordinates by $p(\frac{1}{2}+u,v)=0$.
We refer to this equation as the {\it{canonical equation}} of $\Gamma$. It is immediate to check that
any nontrivial homogeneous part in the corresponding polynomial in $u,v$ has degree with the same parity as $n$;
for instance, on a smooth $4$-fold $\mathcal X$, for any divisor $D = \frac{1}{2}K_{\mathcal X} + \Delta$ the Riemann--Roch formula gives
\begin{equation} \label{RR}
\chi(D)
=\frac{1}{24}\Delta^4 +\frac{1}{48}\bigg(2c_2(\mathcal X)-K_{\mathcal X}^2\bigg)\cdot \Delta^2
+ \frac{1}{384}\bigg(K_{\mathcal X}^2-4c_2(\mathcal X)\bigg)\cdot K_{\mathcal X}^2 + \chi(\mathcal O_{\mathcal X})
\end{equation}
(e.\ g.\ see \cite[p.\ 292]{BLL}).
We thus see that for a polarized $4$-fold, the polynomial $p$ contains only homogeneous parts of degree $4$ and $2$ in $u,v$ plus the constant term: so, if the latter is zero, then $\Gamma$ has a singular point at the origin.

The most significant property of the Hilbert curve of $(\mathcal X,\mathcal L)$ is its sensitivity with respect to fibrations that
suitable adjoint linear systems to $\mathcal L$ may induce on $\mathcal X$ \cite[Theorem 6.1]{BLS}. This makes scrolls (of any type) very
interesting from the point of view of their Hilbert curves. In fact if $(\mathcal X,\mathcal L)$ is a scroll over $\mathcal Y$, with
$\dim \mathcal Y=m$, then
$\Gamma_{(\mathcal X,\mathcal L)}$ consists of $n-m$ parallel lines plus a curve $C$, of degree $m$, and we can consider the following
question (see \cite[Problem 6.6]{BLS}).

\smallskip

\begin{quest}\label{quest1}
Can $C$ itself be regarded as the Hilbert curve of $\mathcal Y$, polarized by some ample $\mathbb Q$-line bundle ?
\end{quest}

\smallskip

For instance, for scrolls over a smooth curve the answer is positive \cite[Remark 4.1]{L}.
This note is mainly concerned with the answer to Question \ref{quest1} for the $4$-scroll $(X,L)$ described below (see also \cite[p.\ 330]{BS}, \cite{BS0}).
Set $X=\mathbb P_Y(\mathcal F)$ and let $p:X \to Y$ be the projection, where $Y$ is
$\mathbb P^3$ blown-up at a point $w$, $\mathcal F = H^{\oplus 2}$, and $H = \sigma^*\mathcal O_{\mathbb P^3}(3) - e$,
$\sigma: Y \to \mathbb P^3$ standing for the blowing-up and $e \cong \mathbb P^2$ for the
exceptional divisor.
We denote by $L$ the tautological line bundle of $\mathcal F$ on $X$. Clearly,
$(X,L)$ is a classical scroll over $Y$ via $p$, while it is an adjunction theoretic scroll over $\mathbb P^3$ via
the map $\pi:=\sigma \circ p:X \to \mathbb P^3$, since
\begin{equation}\label{beginning}
K_X+2L = p^*(K_Y+2H) = \pi^*\bigg(K_{\mathbb P^3}+ 2 \mathcal O_{\mathbb P^3}(3)\bigg)=\pi^* \mathcal O_{\mathbb P^3}(2).
\end{equation}
However, it is not a classical scroll over $\mathbb P^3$, since the fiber $\pi^{-1}(w) = \mathbb P_e(\mathcal F_e)$
is a divisor inside $X$, being isomorphic to $e \times \mathbb P^1$. The following diagram summarizes the above situation
$$
\xymatrix{
X=\mathbb{P}_Y(\mathcal{F})\ar@{->}[d]_{p} \ar@{->}[dr]^{\pi} \\
Y \ar@{->}[r]^{\sigma} & \mathbb{P}^3 \ \ .
}
$$
Before addressing Question \ref{quest1} for $(X,L)$, we need the equation of $\Gamma_{(X,L)}$.

\begin{prop} Let $(X,L)$ be the pair described above. The canonical equation of $\Gamma_{(X,L)}$ in coordinates $(u,v)$, is
\begin{equation}\label{equation}
p_{(X,L)}\bigg(\frac{1}{2}+u,v\bigg) = \frac{1}{3} (2u-v) (u-v) (28u^2 - 38uv + 13v^2 - 1) = 0.
\end{equation}
\end{prop}

Section \ref{first} and Section \ref{second} contain two different proofs of this statement.

\section{First approach} \label{first}

Here, to get the canonical equation of the Hilbert curve $\Gamma_{(X,L)}$ we implement \eqref{RR} with $\mathcal X = X$ and $\Delta=uK_X+vL$.

\bigskip
First of all we recall the Chern--Wu relation:
$$ L^2 - L \cdot p^*c_1(\mathcal F) + p^*c_2(\mathcal F)=0.$$
Since $\mathcal F = H^{\oplus 2}$, it gives $L^2 =  L \cdot p^*(2H) - p^*(H^2)$. Moreover, for any
divisor $D$ on $Y$, we get
$$L \cdot p^*D^3 = D^3,$$
$$L^2 \cdot p^*D^2 = 2 H \cdot D^2,$$
$$L^3 \cdot p^*D = 3H^2 \cdot D,$$
and $L^4 = 4H^3$. Let $h=\sigma^*\mathcal O_{\mathbb P^3}(1)$; then $H = 3h - e$, hence
$H^3= 26$, since $h^3= \big(\mathcal O_{\mathbb P^3}(1)\big)^3=1$,
$e^3=\big(\mathcal O_e(e)\big)^2=1$, and $h \cdot e=0$. Therefore
\begin{equation}\label{L^4}
L^4 = 104.
\end{equation}
Moreover, specializing the above intersections for $D=h$ and $D=e$ respectively, we get
\begin{equation}\label{L_with_h}
L \cdot p^*h^3 = 1, \qquad L^2 \cdot p^*h^2 = 6, \qquad L ^3\cdot p^*h = 27,
\end{equation}
$$L \cdot p^*e^3 = 1, \qquad L^2 \cdot p^*e^2 = -2, \qquad L ^3\cdot p^*e = 3.$$

Now look at the Chern classes of $X$. We have $K_X = -2L + p^*(K_Y + 2H)$,
by the canonical bundle formula for $X$, since $\mathcal F = H^{\oplus 2}$. Moreover,
since $K_Y=-4h+2e$, we get $K_Y+2H = 2h$, hence
$$K_X = -2L + p^*(K_Y+2H) = -2(L - p^*h).$$
Consequently,
$$K_X^2 = 4 (L^2 - 2L \cdot p^*h + p^*h^2),$$
$$K_X^3 = -8(L^3 - 3L^2 \cdot p^*h + 3L \cdot p^*H^2 - p^*h^3),$$
and
\begin{equation}\label{K^4}
K_X^4 = 16(L^4 - 4L^3 \cdot p^*h + 6 L^2 \cdot p^*h^2 - 4L \cdot p^*h^3)= 16 \times 28 = 448.
\end{equation}
Combining these with \eqref{L_with_h}
we can compute the pluridegrees $K_X^i \cdot L^{4-i}$ for $i=1,2,3$:
\begin{equation}\label{pluridegrees}
K_X \cdot L^3 = -154,    \qquad  K_X^2 \cdot L^2 = 224,    \qquad  K_X^3 \cdot L = -320.
\end{equation}
This provides the values of several intersection products in the Riemann--Roch formula,
but many other involve the second Chern class of $X$.
To evaluate it, looking at the $\mathbb P^1$-bundle structure
$p:X \to Y$, we can use the relative tangent sequence
$$0 \to T_{X/Y} \to T_X \to p^*T_Y \to 0$$
and the relative Euler sequence
$$0 \to \mathcal O_X \to p^* \mathcal F^{\vee} \otimes L \to T_{X/Y} \to 0.$$
Combining them, we get the following relation between the Chern polynomials
$$c(T_X;t) = p^* c(T_Y;t)\ c(p^*\mathcal F^{\vee} \otimes L; t),$$
which gives
$$c_2(X) = p^*c_2(Y) + p^*c_1(Y) \cdot c_1(p^*\mathcal F^{\vee} \otimes L) + c_2(p^*\mathcal F^{\vee} \otimes L).$$
Recall that $c_2(Y)= \sigma^* c_2(\mathbb P^3)$ (e.\ g., see \cite[Lemma at p.\ 609]{GH}),
hence $c_2(Y)=6h^2$. Moreover, $c_1(p^*\mathcal F^{\vee} \otimes L)=2(L-p^*H)$ and
$c_2(p^*\mathcal F^{\vee} \otimes L)=(L-p^*H)^2$. So, taking into account the expressions
of $H$ and $K_Y$ in terms of $h$ and $e$, we obtain
$$c_2(X) = L^2 + 2L \cdot p^*(h-e) - 3p^*(3h^2+e^2).$$
This gives
\begin{equation}\label{c_2K^2}
c_2(X)\cdot K_X^2 = 4\bigg(L^4 - 2L^3 \cdot p^*e - 3L^2 \cdot p^*(4h^2+e^2) + 20 L \cdot p^*h^3  \bigg) = 4 \times 52 = 208.
\end{equation}
Moreover,
$$2c_2(X)-K_X^2 = -2\bigg(L^2 - 2 L \cdot p^*(3h-e) + p^*(11h^2+3e^2) \bigg).$$
As a consequence of the above relations we get
$$\bigg(2c_2(X)- K_X^2\bigg)\cdot K_X^2 = -32, \qquad \bigg(2c_2(X)- K_X^2\bigg)\cdot K_X \cdot L = 24, \qquad \bigg(2c_2(X)- K_X^2\bigg)\cdot L^2 = -16.$$
Now we have all ingredients; so, letting $\Delta = uK_X+vL$, \eqref{RR} allows us to express
the canonical equation of the Hilbert curve $\Gamma_{(X,L)}$. First of all,
since $\chi(\mathcal O_X)=1$, from \eqref{K^4} and \eqref{c_2K^2} we get the degree zero term, which is
$$\frac{1}{384}\bigg(K_X^2-4c_2(X)\bigg)\cdot K_X^2 + \chi(\mathcal O_X)=\frac{1}{384}(448-832)+1 = 0.$$
This means that $\Gamma_{(X,L)}$ has a singular point of multiplicity $\geq 2$ at the origin.
Next, since
$$\bigg(2c_2(X)-K_X^2\bigg) \cdot \Delta^2= -32 u^2 + 48uv -16v^2 = -16(2u^2 - 3uv + v^2) = -16(2u-v)(u-v),$$
in view of the previous computations, the homogeneous part of degree 2 is
\begin{eqnarray} \nonumber
\frac{1}{48}\bigg(2c_2(X)-K_X^2\bigg) \cdot \Delta^2 = -\ \frac{1}{3}(2u-v)(u-v). \nonumber
\end{eqnarray}
As to the homogeneous part of degree 4, \eqref{K^4}, \eqref{pluridegrees} and \eqref{L^4} show that
$$\Delta ^4 = (uK_X+vL)^4 = 8\ F(u,v),$$
where
$$F(u,v) = 56u^4 -160u^3v + 168 u^2v^2 - 77 uv^3 + 13 v^4;$$
hence $\frac{1}{24}\Delta^4 = \frac{1}{3}F(u,v)$.
Note that the polynomial $F$ can be rewritten as
$$F(u,v) = 28u^2 (2u^2-3uv+v^2) - v\ G(u,v)=28 u^2(2u-v)(u-v) - v\ G(u,v),$$
where $G(u,v) = 76u^3 -140u^2v + 77uv^2-13v^3$; moreover, it is easy to see that
\begin{eqnarray} \nonumber 
G(u,v) &=& (2u-v) (38u^2 - 51uv + 13 v^2)\\ \nonumber
&=& (2u-v) (u-v) (38u - 13 v).
\end{eqnarray}
Thus
$$F(u,v) = (2u-v)(u-v) \big[28u^2 - v(38u-13v) \big].$$
Therefore, the homogeneous part of degree $4$ is
$$\frac{1}{24}\Delta^4 = \frac{1}{3} (2u-v)(u-v)(28u^2-38uv+13v^2).$$
In conclusion, putting all pieces together and collecting all common factors, we get
\eqref{equation}.

\section{Second approach}\label{second}

In this section, we obtain equation
\eqref{equation} again
with another approach using Algorithm 3 in \cite[Appendix]{LT}.
To do that, it is more convenient to use coordinates $(x,y)=(\frac{1}{2}+u,v)$ in the plane of $\Gamma_{(X,L)}$.
Let $S\subset\mathbb{P}^3$ be a smooth quadric surface not containing the point $w$ and consider the smooth threefold $V:=\pi^{-1}(S)\in |\pi^*\big(\mathcal{O}_{\mathbb{P}^3}(2)\big)|$.
Clearly, $V\cap\pi^{-1}(w)=\emptyset$ and $(V,L_V)$ is a scroll over $S$ via $\pi|_V:V\to S$, with $K_V+2L_V=(\pi|_V)^*\mathcal{O}_{\mathbb{P}^3}(4)_S$.
Note also that
\begin{equation}\label{V}
V\in |K_X+2L|
\end{equation}
in view of \eqref{beginning}. According to the above quoted algorithm, consider the following exact sequence:
$$0\to xK_X+yL+(x-1)V\to x(K_X+V)+yL\to xK_V+yL_V\to 0,$$
which by \eqref{V} can be rewritten as
\begin{equation}\label{exseq}
0 \to (2x-1)K_X+(2(x-1)+y)L \to 2xK_X+(2x+y)L \to xK_V+yL_V\to 0.
\end{equation}
The exact sequence \eqref{exseq} gives the following relation between $p_{(X,L)}$ and $p_{(V,L_V)}$:
\begin{equation}\label{p}
p_{(X,L)}(2x,2x+y)=p_{(X,L)}(2x-1,2x+y-2)+p_{(V,L_V)}(x,y).
\end{equation}
By \cite[Theorem 6.1]{BLS} we know that, in terms of coordinates $(a,b)$, the two polynomials can be written as
$$p_{(X,L)}(a,b)=R_{(X,L)}(a,b)\cdot (2a-b-1)\qquad \mathrm{and} \qquad p_{(V,L_V)}(a,b)=R_{(V,L_V)}(a,b)\cdot (2a-b-1),$$
where $R_{(X,L)}$ and $R_{(V,L_V)}$ are polynomials in $(a,b)$ of degrees $3$ and $2$, respectively.
Thus \eqref{p} becomes
\begin{equation}\label{R}
R_{(X,L)}(2x,2x+y)=R_{(X,L)}(2x-1,2x+y-2)+R_{(V,L_V)}(x,y).
\end{equation}
The goal will be to find the explicit expression of the polynomial
\begin{equation}\label{RX}
R_{(X,L)}(a,b):=A'a^3+B'a^2b+C'ab^2+E'b^3+F'a^2+G'ab+H'b^2+J'a+L'b+M',
\end{equation}
with rational coefficients, because $R_{(V,L_V)}$ is known by \cite[Theorem 4.3]{LT}.
Actually, adapting the notation used there (see also \cite[Example 4.2]{LT}) to our situation, we have
$$S\cong\mathbb{P}^1\times\mathbb{P}^1, \quad  A=\mathcal{O}_{\mathbb{P}^3}(4)_S, \quad K_S=\mathcal{O}_{\mathbb{P}^3}(-2)_S.$$
Hence $\chi (\mathcal{O}_V)=\chi (\mathcal{O}_S)=1$, $K_S\cdot A=-16$ and $A^2=32$. Moreover, from the exact sequence
$$0 \to L-V=L+\pi^*\mathcal{O}_{\mathbb{P}^3}(-2) \to L \to L_V \to 0,$$
we get $\chi (L_V)=\chi (L)-\chi (L+\pi^*\mathcal{O}_{\mathbb{P}^3}(-2))$. Observe that
$$\chi (L)=\chi (\mathcal{F})=2\chi (H)=2h^0(Y,\sigma^*\mathcal{O}_{\mathbb{P}^3}(3)-e)=38$$ and
$$\chi (L+\pi^*\mathcal{O}_{\mathbb{P}^3}(-2))=\chi (\mathcal{F}\otimes\sigma^*\mathcal{O}_{\mathbb{P}^3}(-2))=2\chi(H\otimes\sigma^*\mathcal{O}_{\mathbb{P}^3}(-2))=
2\chi(\sigma^*\mathcal{O}_{\mathbb{P}^3}(1)-e)=6.$$ Therefore, $\chi (L_V)=38-6=32$ and then from \cite[Theorem 4.3]{LT} we deduce that
\begin{equation}\label{RV}
R_{(V,L_V)}(x,y)=-4x^2+12xy-9y^2+4x-6y-1.
\end{equation}
Note that Serre duality on $X$ implies that $p_{(X,L)}(a,b)=p_{(X,L)}(1-a,-b)$, which in turn gives
$$R_{(X,L)}(a,b)=-R_{(X,L)}(1-a,-b).$$
This leads by using MAPLE to the following relations:
$$A'=2J'+4M', \quad B'=-G', \quad C'=-2H', \quad F'=-3J'-6M'.$$
Using these relations, \eqref{RV} and \eqref{RX} with the pairs $(2x,2x+y)$ and $(2x-1,2x+y-2)$ instead of $(a,b)$ to obtain
the terms $R_{(X,L)}(2x,2x+y)$ and $R_{(X,L)}(2x-1,2x+y-2)$, respectively, from \eqref{R}
we deduce the following expressions for four further unknown coefficients:
$$G'=-12E'-30, \quad H'=3E'+\frac{9}{2}, \quad J'=-4E'-2M'-\frac{38}{3}, \quad L'=2E'+M'+\frac{9}{2}.$$
Finally, by computing $p_{(X,L)}$ in $(0,0)$ and $(0,1)$, we get
$$1=p_{(X,L)}(0,0)=(-1)R_{(X,L)}(0,0)=-M',$$
$$38=\chi (L)=p_{(X,L)}(0,1)=(-2)R_{(X,L)}(0,1)=-12E'-14.$$
Hence $M'=-1$ and $E'=-\frac{13}{3}$. By replacing these values in the previous expressions of the coefficients, we deduce the final expression of $R_{(X,L)}$ in terms of
the coordinates $(x,y)$:
$$R_{(X,L)}(x,y)=\frac{28}{3}x^3-22x^2y+17xy^2-\frac{13}{3}y^3-14x^2+22xy-\frac{17}{2}y^2
+\frac{20}{3}x -\frac{31}{6}y-1,$$
which leads to $p_{(X,L)}(\frac{1}{2}+u,v)$ as in \eqref{equation}, keeping in mind that $p_{(X,L)}(x,y) = R_{(X,L)}(x,y) \cdot (2x-y-1)$
and $(x,y)=\big(\frac{1}{2}+u,v\big)$.

\bigskip

\section{A singular property of $\Gamma_{(X,L)}$}

Coming back to Question \ref{quest1},
we highlight an intriguing property of the Hilbert curve of our polarized fourfold $(X,L)$.
As observed, we can regard $(X,L)$ as an adjunction theoretic scroll over $\mathbb P^3$ as well as a classical scroll over $Y$.
Due to \eqref{beginning},
since \cite[Theorem 6.1]{BLS} holds for scrolls of both types, the linear factor $(2u-v)$ in \eqref{equation} was
{\it{a priori}} expected.
The question is whether the residual degree $3$ factor
$$\phi(u,v) = \frac{1}{3} (28u^3-66u^2v+51uv^2-13v^3 - u + v),$$
defining a plane cubic $C$,
is somehow related to the base threefold $\mathbb P^3$ ($Y$, respectively) of our scroll $(X,L)$ for some
polarization. Let us start with $\mathbb P^3$. By (\ref{p1}) with $x=u+\frac{1}{2}$ and $y=v$, we see that for any positive integer $a$ the canonical equation of the Hilbert curve of
the polarized threefold $\big(\mathbb P^3, \mathcal O_{\mathbb P^3}(a)\big)$ is
\begin{equation}\nonumber
\frac{1}{6} \prod_{i=1}^3 (-4u +av+i-2) = 0,
\end{equation}
and the same occurs for any positive $a\in \mathbb Q$.
It is immediate to check that the polynomial on the left hand contains nontrivial homogeneous terms of degree $2$,
contrary to what happens
for $\phi$. Therefore the cubic $C$ of equation $\phi(u,v)=0$ cannot be the Hilbert curve of
$\big(\mathbb P^3, \mathcal O_{\mathbb P^3}(a)\big)$.

\bigskip

This shows that in general for an adjunction theoretic scroll, Question \ref{quest1}
has a negative answer.

\bigskip

Next consider $Y$. Any ample line bundle $M$ on $Y$
can be written as $M = ah - re$ for suitable integers $a$ and $r$. For any divisor $D'=\frac{1}{2}K_Y + \Delta'$ on $Y$,
the Riemann--Roch formula says that
\begin{equation}\label{rr}
\chi(D') = \frac{1}{6}{\Delta'}^3 + \frac{1}{24}(2c_2(Y) - K_Y^2) \cdot \Delta'
\end{equation}
\cite[p.\ 291]{BLL}. Hence, letting $\Delta'= uK_Y + vM$ and
computing all required intersections, \eqref{rr} leads to the canonical equation of the Hilbert curve of
$(Y,M)$, which turns out to be
$$p_{(Y,M)}\bigg(\frac{1}{2} + u,v\bigg) = \frac{1}{6}\bigg[-56u^3+12(4a-r)u^2v-6(2a^2-r^2)uv^2+(a^3-r^3)v^3 +2u -(a-r)v \bigg]=0.$$
This polynomial is proportional to $\phi(u,v)$ if and only if the matrix
\[ \begin{pmatrix}
-56 & 12(4a-r) & -6(2a^2-r^2) & a^3-r^3 & 2 & r-a \\
28 & -66 & 51 & -13 & -1 & 1\\
\end{pmatrix}
\]
has rank 1. An immediate check shows that this happens if and only if $(a,r)=(3,1)$, i.e. for $M=H$.
We thus see that
$$\phi(u,v) = - \ p_{(Y,H)}\bigg(\frac{1}{2}+u,v\bigg).$$
Therefore, the factor $\phi$ defines the Hilbert curve of the base $Y$ of our classical scroll $(X,L)$,
endowed with the average polarization $H = \frac{1}{2} \det \mathcal F$ induced by the ample vector bundle $\mathcal F$.

\bigskip

Moreover, we see that, in the special situation we are dealing with, Question \ref{quest1}
has a positive answer
regarding $(X,L)$ as a classical scroll over $Y$, while this is not the case when we look at it as an
adjunction theoretic scroll over $\mathbb P^3$.

\bigskip

\noindent {\textit{Remark}}. The conclusion concerning $(X,L)$ as a scroll over $\mathbb P^3$ can be obtained
more geometrically, arguing as follows. The Hilbert curve of $\big(\mathbb P^3, \mathcal O_{\mathbb P^3}(a)\big)$
consists of three parallel evenly spaced lines,
while, from the real point of view, the cubic $C$ consists of the line $u-v=0$ plus an ellipse: actually, a straightforward verification
shows that the conic $\gamma$ of equation $28u^2-38uv+13v^2-1 = 0$ is an ellipse
whose axes, determined by the eigenvectors of the matrix
\begin{equation} \label{Ainfty}
A_{\infty} := \begin{pmatrix}
28  & -19\\
-19 & 13 \\
\end{pmatrix}, \quad
\end{equation}
are $3u-2v=0$ and $2u+3v=0$. From another perspective,
removing both linear factors $2u-v$ and $u-v$ from \eqref{equation} one could ask whether the conic
$\gamma$ described by the residual degree 2 polynomial is the Hilbert curve of some polarized
or $\mathbb Q$-polarized surface $(S,\mathcal L)$.
Even in this case the answer is negative. Otherwise, taking into account that the canonical equation of the Hilbert curve
of $(S,\mathcal L)$ is
$$\frac{1}{2} \Bigg(K_S^2 u^2 +2K_S \cdot \mathcal L uv + \mathcal L^2 v^2 + \big(2\chi(\mathcal O_S) - \frac{1}{4}K_S^2 \big) \Bigg) = 0,$$
\eqref{Ainfty} would imply the existence of a nonzero rational number $\rho$ such that
$$K_S^2 \mathcal L^2 - (K_S \cdot \mathcal L)^2 = \rho^2 \det A_{\infty} = 3\ \rho^2 > 0,$$
but this contradicts the Hodge index theorem.

\vspace{1.5cm}

\noindent {\bf Acknowledgements.} The first author is a member of G.N.S.A.G.A. of the Italian INdAM.
He would like to thank the PRIN 2015 Geometry of Algebraic Varieties and the University
of Milano for partial support. During the preparation of this paper, the second author
was partially supported by the National Project Anillo ACT 1415 PIA CONICYT and the
Proyecto VRID N.214.013.039-1.OIN of the University of Concepci\'on. The authors are grateful to the referee 
for useful remarks.

\end{document}